\documentclass[11pt]{article}

\usepackage{latexsym}
\usepackage{amssymb}

\newtheorem{theorem}{Theorem}

\begin{document}
{
\begin{center}
{\Large\bf
On a family of hypergeometric Sobolev orthogonal polynomials on the unit circle.}
\end{center}
\begin{center}
{\bf S.M. Zagorodnyuk}
\end{center}

\section{Introduction.}

The theories of orthogonal polynomials on the real line and on the unit circle have many similarities as well as considerable 
differences~\cite{cit_50000_Gabor_Szego}, \cite{cit_5000_Ismail}, \cite{cit_48000_Simon_1}, \cite{cit_48000_Simon_2}.
For a long time they have been developed side by side by efforts of numerous mathematicians.
The theory of Sobolev orthogonal polynomials is a much more \textit{terra incognita}~\cite{cit_5150_M_X}, \cite{cit_49000_Sri_Ranga}. 
In this theory one can also see that
some ideas come from the real line to the unit circle. Examples of such ideas are adding of Dirac deltas to the classical inner
products and considering of coherent pairs of measures (see, e.g.,~\cite{cit_5000_Castillo}, \cite{cit_5500_GMPC_coherent_pairs} and
references therein).
In the present paper we shall follow the same line: we shall develop the ideas from~\cite{cit_80000_Zagorodnyuk_JAT_2020}  to get some 
new hypergeometric polynomials and study their properties.

Let $\mu$ be a (non-negative) measure on $\mathfrak{B}(\mathbb{T})$ with an infinite support, $\mu(\mathbb{T})=1$. Denote by
$p_n$ orthogonal polynomials on $\mathbb{T}$ with respect to $\mu$ ($\deg p_n = n$, but the positivity of leading coefficients
is not assumed):
\begin{equation}
\label{f1_10}
\int_{\mathbb{T}} p_n(z) \overline{p_m(z)} d\mu = A_n \delta_{n,m},\qquad A_n>0,\ n,m\in\mathbb{Z}_+.
\end{equation}
Fix an arbitrary positive integer $\rho$. Consider the following differential equation:
\begin{equation}
\label{f1_20}
( e^{-x} y(x) )^{ (\rho) } = e^{-x} p_n(x),\qquad  n\in\mathbb{Z}_+.
\end{equation}
Expanding the derivative by the Leibniz formula and canceling $e^{-x}$ we get
\begin{equation}
\label{f1_25}
\sum_{k=0}^\rho (-1)^{\rho -k} \left( \begin{array}{cc} \rho \\
k \end{array} \right) 
y^{ (k) }(x)
= p_n(x),\qquad  n\in\mathbb{Z}_+.
\end{equation}

\noindent
\textbf{Condition A.}
\textit{Suppose that for each $n\in\mathbb{Z}_+$, there exists a $n$-th degree polynomial solution $y=y_n(x)$ of~(\ref{f1_25})}.

If Condition~A is satisfied, then $y_n$ are Sobolev orthogonal polynomials on $\mathbb{T}$:
\begin{equation}
\label{f1_30}
\int_{\mathbb{T}} \left( y_n(z), y_n'(z),..., y_n^{(\rho)}(z) \right) M \overline{
\left( \begin{array}{cccc} y_m(z) \\
y_m'(z) \\
\vdots \\
y_m^{(\rho)}(z) \end{array} \right)
} 
d\mu = A_n \delta_{n,m},\qquad n,m\in\mathbb{Z}_+,
\end{equation}
where
\begin{equation}
\label{f1_35}
M =
\left( 
(-1)^{l+j}
\left( \begin{array}{cc} \rho \\
l \end{array} \right)
\left( \begin{array}{cc} \rho \\
j \end{array} \right)
\right)_{l,j=0}^\rho.
\end{equation}

In this paper we shall only consider the following case: $p_n(x) = x^n$, and $\mu$ being the normalized arc length measure on $\mathbb{T}$.
Let us briefly describe the content of the paper.
Equating coefficients of the same powers on the both sides of equation~(\ref{f1_25}) one obtains a linear system of equations for 
the coefficients of an unknown polynomial $y(x)$ (the same idea was used in~\cite{cit_5_Azad}). However, for large values of $\rho$
it is not easy to get a convenient expression for solutions, without huge determinants or recurrences. In this case equation~(\ref{f1_20})
turned out to be useful. It gives a possibility to express $y_n(x)$ for $\rho=1$ in terms of the incomplete gamma function.
A step-by-step analysis for $\rho = 1,2,...,$ allows to obtain an explicit representation of $y_n(x)$. Explicit representations, 
differential equations and orthogonality relations for $y_n$ will be given by Theorem~1. Using Fasenmeier's method (\cite{cit_5150_Rainville}) for
the reversed polynomials $y_n^*(x) = x^n y_n\left( \frac{1}{x} \right)$, we shall derive recurrence relations for $y_n(x)$ (Theorem~2) as well.

\noindent
{\bf Notations. }
As usual, we denote by $\mathbb{R}, \mathbb{C}, \mathbb{N}, \mathbb{Z}, \mathbb{Z}_+$,
the sets of real numbers, complex numbers, positive integers, integers and non-negative integers,
respectively. 
Set $\mathbb{T} := \{ z\in\mathbb{C}:\ |z|=1 \}$. By $\mathfrak{B}(\mathbb{T})$ we mean the set of all Borel subsets of $\mathbb{T}$.
By $\mathbb{P}$ we denote the set of all polynomials with complex coefficients.
For a complex number $c$ we denote
$(c)_0 = 1$, $(c)_1=c$, $(c)_k = c(c+1)...(c+k-1)$, $k\in\mathbb{N}$ (\textit{the shifted factorial or Pochhammer symbol}).
The generalized hypergeometric function is denoted by
$$ {}_m F_n(a_1,...,a_m; b_1,...,b_n;x) = \sum_{k=0}^\infty \frac{(a_1)_k ... (a_m)_k}{(b_1)_k ... (b_n)_k} \frac{x^k}{k!}, $$
where $m,n\in\mathbb{N}$, $a_j,b_l\in\mathbb{C}$.

\section{Some Sobolev orthogonal polynomials on $\mathbb{T}$.}

As it was stated in the Introduction, in what follows we shall consider the following case:
$p_n(x) = x^n$, and $\mu=\mu_0$ being the (probability) normalized arc length measure on $\mathbb{T}$. Rewrite equations~(\ref{f1_20}),(\ref{f1_25})
for this case:
\begin{equation}
\label{f2_20}
( e^{-x} y_n(x) )^{ (\rho) } = e^{-x} x^n,\qquad  n\in\mathbb{Z}_+;
\end{equation}
\begin{equation}
\label{f2_25}
\sum_{k=0}^\rho (-1)^{\rho -k} \left( \begin{array}{cc} \rho \\
k \end{array} \right) 
y_n^{ (k) }(x)
= x^n,\qquad  n\in\mathbb{Z}_+.
\end{equation}
We start with the case $\rho=1$. In this case equation~(\ref{f2_25}) has the following form:
\begin{equation}
\label{f2_27}
y_n'(x) - y_n(x) = x^n,\qquad  n\in\mathbb{Z}_+.
\end{equation}
Fix an arbitrary $n\in\mathbb{Z}_+$. We shall seek for a solution of the required form:
\begin{equation}
\label{f2_29}
y_n(x) = \sum_{k=0}^n \mu_{n,k} x^k,\qquad \mu_{n,k}\in\mathbb{C}.
\end{equation}
Substitute for $y_n$ into~(\ref{f2_27}) to get
$$ \sum_{k=0}^{n-1} 
\{
(k+1) \mu_{n,k+1} - \mu_{n,k}
\}
x^k
-
\mu_{n,n} x^n = x^n. $$
Comparing the coefficients of the same powers on the both sides we obtain that
\begin{equation}
\label{f2_31}
\mu_{n,n} = -1,\quad \mu_{n,k} = (k+1) \mu_{n,k+1},\quad k=n-1,n-2,...,0.
\end{equation}
It can be verified by the induction argument that
$$ \mu_{n,j} = - (n)_{n-j} = - \frac{n!}{j!},\qquad j=0,1,...,n. $$
Thus,
\begin{equation}
\label{f2_35}
y_n(x) = - n! \sum_{k=0}^n \frac{ x^k }{ k! },
\end{equation}
is a solution of~(\ref{f2_27}). In the case $\rho > 1$, it is not easy to solve the corresponding recurrence relation for the coefficients
and we shall proceed in another way.

Observe that
\begin{equation}
\label{f2_37}
y_n(x) = - e^x \Gamma(n+1,x),
\end{equation}
where
$$ \Gamma(\alpha,x) = \int_x^\infty e^{-t} t^{\alpha - 1} dt,\qquad \alpha > 0, $$
is \textit{the complementary incomplete gamma function} (\cite{cit_3_Andrews_book}).
In fact, integrating~(\ref{f2_20}) (with $\rho=1$) from $a$ to $b$ ($a,b\in\mathbb{R}$) we get
$$ e^{-b} y_n(b) - e^{-a} y_n(a) = \int_a^b e^{-x} x^n dx. $$
Passing to the limit as $b\rightarrow +\infty$ we get
\begin{equation}
\label{f2_38}
y_n(a) = - e^a \int_a^\infty e^{-x} x^n dx,
\end{equation}
and relation~(\ref{f2_37}) follows.

Suppose that we have constructed a polynomial solution (of the required form) $y_n(\rho;x) = y_n(x)$ of equation~(\ref{f2_20})
for some positive integer $\rho$. Let us show how to get a polynomial solution $y_n(\rho+1;x)$ of equation~(\ref{f2_20}) with $\rho+1$.
Notice that we do not state the uniqueness of such solutions for $\rho>2$.  
We shall need the following auxiliary equation:
\begin{equation}
\label{f2_42}
( e^{-x} y_n(\rho+1;x) )' = e^{-x} y_n(\rho;x),\qquad  n\in\mathbb{Z}_+,
\end{equation}
with an unknown $y_n(\rho+1;x)$.
Equation~(\ref{f2_42}) has a unique $n$-th degree polynomial solution. This can be verified comparing the coefficients of polynomials,
in the same manner as for equation~(\ref{f2_27}).
It is not easy to solve the corresponding recurrence relation in this case, but the existence and the uniqueness of a $n$-th degree polynomial solution
is obvious. 

Integrating relation~(\ref{f2_42}) from $t$ to $b$ we get
$$ e^{-b} y_n(\rho+1;b) - e^{-t} y_n(\rho+1;t) = \int_t^b e^{-x} y_n(\rho;x) dx. $$
Passing to the limit as $b\rightarrow +\infty$ we get
\begin{equation}
\label{f2_44}
y_n(\rho+1;t) = - e^t \int_t^\infty e^{-x} y_n(\rho;x) dx. 
\end{equation}
By~(\ref{f2_20}),(\ref{f2_42}) we may write:
$$ e^{-x} x^n = ( e^{-x} y_n(\rho;x) )^{ (\rho) } = ( e^{-x} y_n(\rho+1;x) )^{ (\rho+1) }. $$
Therefore $y_n(\rho+1;\cdot)$ given by~(\ref{f2_44}) is a required polynomial solution of~(\ref{f2_20}) for $\rho+1$.

Equation~(\ref{f2_44}) shows how to construct polynomial solutions step by step for $\rho = 1,2,...$. 
However, we are interested to get an explicit representation for every $y_n(\rho;x)$.
Let
\begin{equation}
\label{f2_46}
y_n(\rho;x) = \sum_{j=0}^n d_j(\rho) \frac{x^j}{j!},\qquad n\in\mathbb{Z}_+,\ \rho\in\mathbb{N},
\end{equation}
with some unknown complex numbers $d_j(\rho)$.
By~(\ref{f2_44}),(\ref{f2_38}),(\ref{f2_35}) we may write
$$ y_n(\rho+1;t) = - e^t \sum_{j=0}^n d_j(\rho) \frac{1}{j!} \int_t^\infty e^{-x} x^j dx = 
\sum_{j=0}^n d_j(\rho) \frac{1}{j!} y^j(1;t) = $$
\begin{equation}
\label{f2_50}
= -\sum_{j=0}^n \sum_{k=0}^j d_j(\rho) \frac{x^k}{k!},\qquad n\in\mathbb{Z}_+,\ \rho\in\mathbb{N}.
\end{equation}
Changing the order of summation in~(\ref{f2_50}) we write:
$$ y_n(\rho+1;t) = -\sum_{k=0}^n \sum_{j=k}^n d_j(\rho) \frac{x^k}{k!}. $$
Therefore
\begin{equation}
\label{f2_52}
d_k(\rho+1) = -\sum_{j=k}^n d_j(\rho),\qquad j=0,1,...,n;\ \rho\in\mathbb{N}.
\end{equation}
Relation~(\ref{f2_52}) can be written in a matrix form for the vectors of coefficients
$\vec d(\rho) := (d_0(\rho),...,d_n(\rho))^T$, and a $(n+1)\times(n+1)$ upper-diagonal Toeplitz matrix $T$,
having all nonzero elements equal to $1$:
\begin{equation}
\label{f2_54}
\vec d(\rho+1) = - T  \vec d(\rho),\qquad \rho\in\mathbb{N}.
\end{equation}
Therefore
\begin{equation}
\label{f2_56}
\vec d(\rho) = (-1)^\rho T^\rho (0,...,0,1)^T,\qquad \rho\in\mathbb{N}.
\end{equation}
Applying the Riesz calculus for evaluating $T^\rho$, one obtains the following solution:
\begin{equation}
\label{f2_58}
d_k(\rho) = (-1)^\rho \left( \begin{array}{cc} n-k+\rho-1 \\
n-k \end{array} \right),\qquad k=0,1,...,n;\ n\in\mathbb{Z}_+,\ \rho\in\mathbb{N}.
\end{equation}
We shall omit the details of calculating the resolvent $(T-\lambda E)^{-1}$. We only notice that it was convenient to subtract the subsequent rows
when solving the linear system of equations $(T-\lambda E) f = (0,...,0,1)^T$.
It can be directly verified that the resulting expression~(\ref{f2_58}) obeys~(\ref{f2_52}), by using the Pascal identity and the induction
argument. 

Thus, we have obtained the following representation for $y_n$:
\begin{equation}
\label{f2_60}
y_n(\rho;x) = (-1)^\rho \sum_{j=0}^n \left( \begin{array}{cc} n-k+\rho-1 \\
n-k \end{array} \right)
\frac{x^j}{j!},\qquad n\in\mathbb{Z}_+,\ \rho\in\mathbb{N}.
\end{equation}

\begin{theorem}
\label{t2_1}
Let $y_n(\rho;x)$ be polynomials given by relation~(\ref{f2_60}) ($\rho\in\mathbb{N}$, $n\in\mathbb{Z}_+$). They have the following properties:

\begin{itemize}
\item[(a)] Polynomials $y_n(\rho;x)$ admit the following representation:
\begin{equation}
\label{f2_62}
y_n(\rho;x) = \frac{ (-1)^\rho }{n!} x^n {}_2 F_0 \left(-n,\rho;-;-\frac{1}{x}\right),\quad n\in\mathbb{Z}_+,\ \rho\in\mathbb{N};\ 
x\in\mathbb{C}\backslash\{ 0 \}.
\end{equation}

\item[(b)] Polynomials $y(x) = y_n(\rho;x)$ satisfy the following differential equation:
\begin{equation}
\label{f2_64}
x 
\sum_{k=0}^\rho (-1)^{\rho -k} \left( \begin{array}{cc} \rho \\
k \end{array} \right) 
y^{ (k+1) }(x)
-
n
\sum_{k=0}^\rho (-1)^{\rho -k} \left( \begin{array}{cc} \rho \\
k \end{array} \right) 
y^{ (k) }(x)
= 0.
\end{equation}

\item[(c)] Polynomials $y(x) = y_n(\rho;x)$ obey the following differential equation:
\begin{equation}
\label{f2_65}
x y''(x) - (x+\rho-1) y'(x) - n
\left[
y'(x) - y(x)
\right]
= 0.
\end{equation}

\item[(d)] Polynomials $y_n(x) = y_n(\rho;x)$ are Sobolev
orthogonal polynomials on $\mathbb{T}$:
\begin{equation}
\label{f2_67}
\int_{\mathbb{T}} \left( y_n(z), y_n'(z),..., y_n^{(\rho)}(z) \right) M \overline{
\left( \begin{array}{cccc} y_m(z) \\
y_m'(z) \\
\vdots \\
y_m^{(\rho)}(z) \end{array} \right)
} 
d\mu_0 = \delta_{n,m},\qquad n,m\in\mathbb{Z}_+,
\end{equation}
where $M$ is given by~(\ref{f1_35}).
\end{itemize}

\end{theorem}
\textbf{Proof.} $(a)$: It is readily checked that the reversed polynomial for $y_n$ is given by
$$ y_n^*(\rho;x) = 
\frac{ (-1)^\rho }{n!} x^n {}_2 F_0 \left(-n,\rho;-;-x\right),\qquad n\in\mathbb{Z}_+,\ \rho\in\mathbb{N}, $$
and relation~(\ref{f2_62}) follows.

\noindent
$(b)$: Substitute for $x^n$ from~(\ref{f2_25}) into the following identity:
$$ x (x^n)' = n x^n. $$

\noindent
$(c)$: Hypergeometric polynomials 
\begin{equation}
\label{f2_68}
u = u_n(z) := F_0 \left(-n,\rho;-;z\right),\qquad n\in\mathbb{Z}_+,\ \rho\in\mathbb{N},
\end{equation}
satisfy the following
differential equation:
\begin{equation}
\label{f2_70}
z(-n+\theta)(\rho+\theta) u - \theta u = 0,
\end{equation}
where $\theta = z\frac{d}{dz}$. The differential equation for the generalized hypergeometric function ${}_p F_q$ is usually written
when $p,q\geq 1$. However, the arguments in~\cite[p. 75]{cit_5150_Rainville} can be applied in the case $q=0$ as well.
Then for $z\not=0$ we may write
\begin{equation}
\label{f2_72}
z^2 u''(z) +(\rho+1)z u'(z)- n(zu'(z)+\rho u(z)) - u'(z) = 0.
\end{equation}
Observe that
\begin{equation}
\label{f2_73}
u_n(z) = \frac{n!}{(-1)^{n+\rho}} z^n y_n(\rho;-\frac{1}{z}). 
\end{equation}

Calculating the derivatives $u_n',u_n''$ and inserting them into relation~(\ref{f2_72}), after some algebraic
simplifications, we get relation~(\ref{f2_65}).

\noindent
$(d)$: This follows from our motivation and relation~(\ref{f1_30}) in the Introduction.
$\Box$

In order to obtain a recurrence relation for polynomials $y_n(\rho;x)$ we shall apply Fasenmeier's method~(\cite{cit_5150_Rainville})
to hypergeometric polynomials $u_n(z)$ from~(\ref{f2_68}). In the following considerations, we shall admit for $\rho$
not only positive integer values but $\rho>0$ as well. 
We shall express $u_n,u_{n-1},u_{n-2},zu_n(z),zu_{n-1}(z)$, using $u_{n+1}(z)$.
Choose and fix an arbitrary integer $n$ greater or equal to $2$.
We may write
$$ u_{n+1}(z) = \sum_{k=0}^\infty (-n-1)_k (\rho)_k \frac{z^k}{k!} = \sum_{k=0}^\infty \varepsilon_{n+1}(k), $$
where $\varepsilon_{n+1}(k) = \varepsilon_{n+1}(z;\rho;k) := (-n-1)_k (\rho)_k \frac{z^k}{k!}$.
Using
$$ (-n)_k = (-n-1)_k \frac{ (n+1-k) }{ (n+1) },\qquad k\in\mathbb{Z}_+, $$
$$ (-n+1)_k = (-n-1)_k \frac{ (n+1-k)(n-k) }{ (n+1)n },\qquad k\in\mathbb{Z}_+, $$
and similar relations we obtain that
\begin{equation}
\label{f2_74}
u_n(z) = \sum_{k=0}^\infty \varepsilon_{n+1}(k) \frac{ (n+1-k) }{ (n+1) },
\end{equation}
\begin{equation}
\label{f2_76}
u_{n-1}(z) = \sum_{k=0}^\infty \varepsilon_{n+1}(k) \frac{ (n+1-k)(n-k) }{ (n+1)n },
\end{equation}
\begin{equation}
\label{f2_78}
u_{n-2}(z) = \sum_{k=0}^\infty \varepsilon_{n+1}(k) \frac{ (n+1-k)(n-k)(n-1-k) }{ (n+1)n(n-1) },
\end{equation}
\begin{equation}
\label{f2_80}
z u_n(z) = \sum_{k=0}^\infty \varepsilon_{n+1}(k) \frac{ (-k) }{ (n+1)(\rho+k-1) },\quad \rho\not=1,
\end{equation}
\begin{equation}
\label{f2_82}
z u_{n-1}(z) = \sum_{k=0}^\infty \varepsilon_{n+1}(k) \frac{ (-k)(n+1-k) }{ n(n+1)(\rho+k-1) },\quad \rho\not=1.
\end{equation}
We now assume that $\rho\not=1$.
Consider the following expression $R_n(z)$:
$$ R_n(z) := \varphi_1 u_{n-1}(z) + \varphi_2 u_{n}(z) + \varphi_3 u_{n+1}(z) + 
\varphi_4 z u_{n}(z) + $$
\begin{equation}
\label{f2_84}
+ \varphi_5 z u_{n-1}(z) + \varphi_6 u_{n-2}(z),\qquad \varphi_k\in\mathbb{C}. 
\end{equation}
We intend to choose parameters $\varphi_k$ (depending on the chosen $n$) in such a way that $R_n(z)= 0$, $\forall z\in\mathbb{C}$ .
Substitute above expressions for $u_{n-2},u_{n-1},u_n,zu_n,zu_{n-1}$ into~(\ref{f2_84}) to get
$$ R_n(z) = \sum_{k=0}^\infty \varepsilon_{n+1}(k) \frac{1}{ (n-1)n(n+1)(\rho+k-1) } I_{n,k}, $$
where
$$ I_{n,k} = 
\varphi_1 (n-k)(n+1-k) (n-1)(\rho+k-1) + 
\varphi_2 (n+1-k) (n-1)n(\rho+k-1) +
$$
$$ + \varphi_3 (n-1)n(n+1)(\rho+k-1) + \varphi_4 (-1) k (n-1)n + \varphi_5 (n+1-k) (-1)k (n-1) +$$
\begin{equation}
\label{f2_86}
+ \varphi_6 (n+1-k)(n-k)(n-1-k) (\rho+k-1).
\end{equation}
Observe that $I_{n,k}$ is a polynomial of degree $\leq 4$. Therefore we may check that $I_{n,k}=0$
for some distinct five values of $k$ to get $R_n(z)\equiv 0$. This is a crucial point in Fasenmeier's
method.

We choose $k=-\rho+1; n+1; n; n-1; 0$. After some obvious simplifications we get the following five equations:
\begin{equation}
\label{f2_88}
\varphi_5 = -\frac{ \varphi_4 }{ n+\rho },
\end{equation}
\begin{equation}
\label{f2_90}
\varphi_3 = \frac{ \varphi_4 }{ n+\rho },
\end{equation}
\begin{equation}
\label{f2_92}
\varphi_2 (\rho+n-1) + \varphi_3 (n+1)(\rho+n-1) + \varphi_4 (-1)n - \varphi_5 = 0,
\end{equation}
$$ \varphi_1 2(\rho+n-2) + \varphi_2 2n (\rho+n-2) + \varphi_3 n(n+1) (\rho+n-2) - $$
\begin{equation}
\label{f2_94}
- \varphi_4 (n-1)n - \varphi_5 2(n-1) = 0,
\end{equation}
\begin{equation}
\label{f2_96}
\varphi_1 + \varphi_2 + \varphi_3 + \varphi_6 = 0.
\end{equation}
Set $\varphi_4 = n+\rho$. Then
$$ \varphi_5 = -1,\quad \varphi_3 = 1. $$
By~(\ref{f2_92}) we get
$$ \varphi_2 = \frac{ n-\rho }{ n+\rho-1 }. $$
By~(\ref{f2_94}) we obtain that
$$ \varphi_1 = -\frac{ n(n+1) }{ 2 } -\frac{ (n-\rho)n }{ n+\rho-1 } - 
\frac{ n-1 }{ n+\rho-2 } + \frac{ (n+\rho) (n-1)n }{ 2(n+\rho-2) }. $$
Finally, by~(\ref{f2_96}) we conclude that
$$  \varphi_6 = \frac{ n(n+1) }{ 2 } + \frac{ (n-\rho)(n-1) }{ n+\rho-1 } + 
\frac{ n-1 }{ n+\rho-2 } -\frac{ (n+\rho) (n-1)n }{ 2(n+\rho-2) } - 1. $$
Consequently, polynomials $u_n(z)$ satisfy the following relation:
$$ \left(
-\frac{ n(n+1) }{ 2 } -\frac{ (n-\rho)n }{ n+\rho-1 } - 
\frac{ n-1 }{ n+\rho-2 } + \frac{ (n+\rho) (n-1)n }{ 2(n+\rho-2) }
\right)
u_{n-1}(z) + $$
$$ + \frac{ n-\rho }{ n+\rho-1 }
u_{n}(z) + u_{n+1}(z) + 
(n+\rho)
z u_{n}(z) - z u_{n-1}(z) + $$
$$ +
\left(
\frac{ n(n+1) }{ 2 } + \frac{ (n-\rho)(n-1) }{ n+\rho-1 } + 
\frac{ n-1 }{ n+\rho-2 } -\frac{ (n+\rho) (n-1)n }{ 2(n+\rho-2) } - 1
\right) * $$
\begin{equation}
\label{f2_98}
* u_{n-2}(z),\qquad n=2,3,...; \rho>0,\rho\not=1.
\end{equation}

\begin{theorem}
\label{t2_2}
Let $y_n(x) = y_n(\rho;x)$ be polynomials from relation~(\ref{f2_60}) with $\rho\in\mathbb{N}\backslash\{ 1\}$. 
They satisfy the following recurrence relation:
$$ \left(
-\frac{ n(n+1) }{ 2 } -\frac{ (n-\rho)n }{ n+\rho-1 } - 
\frac{ n-1 }{ n+\rho-2 } + \frac{ (n+\rho) (n-1)n }{ 2(n+\rho-2) }
\right)
x^2 * $$
$$ * (n-1) y_{n-1}(x) +
\frac{ n-\rho }{ n+\rho-1 }
x (n-1)n y_n(x) +
(n-1)n(n+1) y_{n+1}(x) + $$
$$ + 
\left(
\frac{ n(n+1) }{ 2 } + \frac{ (n-\rho)(n-1) }{ n+\rho-1 } + 
\frac{ n-1 }{ n+\rho-2 } -\frac{ (n+\rho) (n-1)n }{ 2(n+\rho-2) } - 1
\right) x^3 * $$
$$ * y_{n-2}(x) -
(n+\rho) (n-1)n y_n(x) +
x (n-1) y_{n-1}(x) = 0, $$
\begin{equation}
\label{f2_100}
n=2,3,....
\end{equation}

\end{theorem}
\textbf{Proof.} Use relations~(\ref{f2_73}) and (\ref{f2_98}).
$\Box$

Observe that $y_1(\rho;x) = (-1)^\rho (x+\rho)$. Thus $y_1(1;x)$ has its root on the unit circle, while
the roots of $y_1(\rho;x)$, for $\rho>1$, are outside $\mathbb{T}$. Consequently, polynomials $y_n(\rho;x)$
are not orthogonal on the unit circle with respect to a scalar measure.

\begin{center}
{\large\bf 
On a family of hypergeometric Sobolev orthogonal polynomials on the unit circle.}
\end{center}
\begin{center}
{\bf S.M. Zagorodnyuk}
\end{center}

In this paper we study the following family of hypergeometric polynomials:
$y_n(x) = \frac{ (-1)^\rho }{ n! } x^n {}_2 F_0(-n,\rho;-;-\frac{1}{x})$, depending on a parameter $\rho\in\mathbb{N}$.
Differential equations of orders $\rho+1$ and $2$ for these polynomials are given.
A recurrence relation for $y_n$ is derived as well.
Polynomials $y_n$ are Sobolev orthogonal polynomials on the unit circle
with an explicit matrix measure.

\vspace{1.5cm}

V. N. Karazin Kharkiv National University \newline\indent
School of Mathematics and Computer Sciences \newline\indent
Department of Higher Mathematics and Informatics \newline\indent
Svobody Square 4, 61022, Kharkiv, Ukraine

Sergey.M.Zagorodnyuk@gmail.com; Sergey.M.Zagorodnyuk@univer.kharkov.ua

}

\begin{thebibliography}{99}

\bibitem{cit_3_Andrews_book}
Andrews, Larry C. Special functions of mathematics for engineers. Reprint of the 1992 second edition. 
SPIE Optical Engineering Press, Bellingham, WA; Oxford University Press, Oxford, 1998. {\rm xx}+480 pp.

\bibitem{cit_5_Azad}
Azad H., Laradji A., Mustafa M. T. Polynomial solutions of differential equations. Adv. Difference Equ. 2011:58 (2011), 12 pp.


\bibitem{cit_5000_Castillo}
Castillo, Kenier. A new approach to relative asymptotic behavior for discrete Sobolev-type orthogonal polynomials 
on the unit circle. Appl. Math. Lett. 25 (2012), no. 6, 1000--1004.

\bibitem{cit_5500_GMPC_coherent_pairs}
Garza, Luis; Marcell\'an, Francisco; Pinz\'on-Cort\'es, Natalia C. $(1,1)$-coherent pairs on the unit circle. Abstr. Appl. Anal. 
2013, Art. ID 307974, 8 pp.

\bibitem{cit_5000_Ismail}
Ismail, Mourad E. H. Classical and quantum orthogonal polynomials in one variable. With two chapters by Walter Van Assche. 
With a foreword by Richard A. Askey. Encyclopedia of Mathematics and its Applications, 98. Cambridge University Press, Cambridge, 2005. xviii+706 pp.

\bibitem{cit_5150_M_X}
Marcell\'an, Francisco; Xu, Yuan. On Sobolev orthogonal polynomials. Expo. Math. 33 (2015), no. 3, 308--352.

\bibitem{cit_5150_Rainville}
Rainville, Earl D. Special functions. Reprint of 1960 first edition. Chelsea Publishing Co., Bronx, N.Y., 1971. {\rm xii}+365 pp.


\bibitem{cit_48000_Simon_1}
Simon, Barry. Orthogonal polynomials on the unit circle. Part 1. Classical theory. 
American Mathematical Society Colloquium Publications, 54, Part 1. American Mathematical Society, Providence, RI, 
2005. xxvi+466 pp.

\bibitem{cit_48000_Simon_2}
Simon, Barry. Orthogonal polynomials on the unit circle. Part 2. Spectral theory. American Mathematical Society Colloquium Publications, 54, 
Part 2. American Mathematical Society, Providence, RI, 2005. pp. i--xxii and 467--1044.

\bibitem{cit_49000_Sri_Ranga}
Sri Ranga, A. Orthogonal polynomials with respect to a family of Sobolev inner products on the unit circle. 
Proc. Amer. Math. Soc. 144 (2016), no. 3, 1129--1143.

\bibitem{cit_50000_Gabor_Szego}
Szeg\"o, G\'abor. Orthogonal polynomials. Fourth edition. 
American Mathematical Society, Colloquium Publications, Vol. XXIII. American Mathematical Society, Providence, R.I., 1975. xiii+432 pp.

\bibitem{cit_80000_Zagorodnyuk_JAT_2020}
Zagorodnyuk, Sergey M. On some classical type Sobolev orthogonal polynomials. J. Approx. Theory 250 (2020), 105337, 14 pp.


\end{thebibliography}
\end{document}